\documentclass[12pt]{amsart}
\usepackage{amssymb}
\usepackage{amsbsy}
\usepackage{amscd}
\usepackage{mathrsfs}
\usepackage{nicefrac}

\usepackage{a4wide}
\usepackage[pdftex]{pict2e}
\usepackage{curve2e}
\makeatletter
\renewcommand{\thesubsection}{(\@roman\c@subsection)}
\makeatother
\newcounter{number}
\usepackage{refcount}

\usepackage{verbatim}
\usepackage{version}
\usepackage{color}
\definecolor{deepjunglegreen}{rgb}{0.0, 0.29, 0.29}
\newenvironment{NB}{
\color{red}{\bf NB}. \footnotesize
}{}
\excludeversion{NB}
\excludeversion{NB2}
\newtheorem{Theorem}{Theorem}

\theoremstyle{definition}

\theoremstyle{remark}

\numberwithin{equation}{section}
\newcommand{\subsecref}[1]{\S\ref{#1}}

\newcommand{\CC}{{\mathbb C}}

\newcommand{\RR}{{\mathbb R}}
\newcommand{\SU}{\operatorname{\rm SU}}
\newcommand{\GL}{\operatorname{GL}}

\newcommand{\U}{\operatorname{\rm U}}

\newcommand{\SO}{\operatorname{\rm SO}}
\newcommand{\grpSp}{\operatorname{\rm Sp}}

\newcommand{\Hom}{\operatorname{Hom}}

\renewcommand{\MR}[1]{}

\newcommand{\vin}[1]{\operatorname{i}(#1)} % incoming vertex
\newcommand{\vout}[1]{\operatorname{o}(#1)} % outgoing vertex

\newcommand{\tslabar}{\mathbin{
\setbox0=\hbox{/\!\!/\!\!/}\rule[0.4\ht0]{\wd0}{.3\dp0}\kern-\wd0\box0}}

\newcommand{\bv}{\mathbf v}
\newcommand{\bw}{\mathbf w}

\usepackage{url}

\usepackage[%dvipdfmx,
bookmarks=true,
citecolor=deepjunglegreen,
colorlinks=true,
]{hyperref}

\begin{document}

\title{Instantons on ALE spaces for classical groups}
\author[H.~Nakajima]{Hiraku Nakajima}
\address{Research Institute for Mathematical Sciences,
Kyoto University, Kyoto 606-8502,
Japan}
\email{nakajima@kurims.kyoto-u.ac.jp}

\begin{abstract}
    We give an ADHM type description of instantons on ALE spaces for
    classical groups as an extension of the description in \cite{KN}
    for unitary groups.
\end{abstract}
\maketitle

\setcounter{tocdepth}{2}
%\tableofcontents

When a gauge group is a unitary group, instantons on ALE spaces have a
description in terms of quiver representations \cite{KN}. It is a
modification of ADHM description \cite{ADHM} of instantons on
$S^4$. (It followed more closely the presentation in \cite{MR1079726}.)
The latter has a version for $\SO/\grpSp$ gauge groups, as an
$\SO/\grpSp$ instanton can be considered as a unitary instanton $A$
together with an involutive isomorphism between the dual instanton and
the original instanton $A^*\cong A$. More precisely we take ADHM
description for both $A$ and $A^*$, and assign an isomorphism between
two descriptions. This was explained already in \cite{ADHM}, and has
been well-known in gauge theory context. See \cite{MR857374} for
instantons on $\overline{\CC P}^2$. It is discussed e.g.\ in
\cite{MR3508922} in the presentation in \cite{MR1079726}.

When \cite{KN} was written, description of $\SO/\grpSp$ instantons on
ALE spaces was \emph{not} known, as the description of the dual
instanton $A^*$ was not given. More precisely, the description in
\cite{KN} involves the parameter $\zeta$ for the level of the
hyperk\"ahler moment map, and we are forced to change $\zeta$ to
$-\zeta$ when we describe $A^*$.
\begin{NB}
Also a diagram automorphism is
involved, which is easy to handle.
\end{NB}%
Thus we need an isomorphism from description for $\zeta$ to one for
$-\zeta$. It was found later as the reflection functor $\mathcal F$
corresponding to the longest element $w_0$ in the finite Weyl group
\cite[\S9]{Na-reflect}.

But the ADHM description of $\SO/\grpSp$ instantons on ALE spaces was
\emph{not} mentioned in \cite[\S9]{Na-reflect}, as the motivation
there was different. We do this job in this short note. 
The author is motivated to do it by a recent preprint \cite{QuiverSym}.
No new input other than \cite{KN,Na-reflect} is necessary, so we just
state the result without a proof. It is also good to look at
\cite[App.~A.4]{2015arXiv150303676N} where ADHM description of
$\SO/\grpSp$ instantons on $\RR^4/\Gamma$ is explained.
It can be considered as a degenerate case of the discussion below
where the reflection functor $\mathcal F$ becomes the identity.

\subsection{Involution on an affine Dynkin diagram}\label{subsec:diagram}

Let $\Gamma$ be a finite subgroup of $\SU(2)$.
Let us define an involution $*$ on the affine Dynkin diagram by
$\rho_i^*\cong\rho_{i^*}$, where $\rho_i$ is an irreducible
representation of $\Gamma$ corresponding to a vertex $i$ via the McKay
correspondence and $\rho_i^*$ is the dual representation of $\rho_i$.
It fixes the trivial representation $\rho_0$, and hence induces a
diagram involution on the finite Dynkin diagram. It is the same
diagram involution given by the longest element $w_0$ of the Weyl
group as $-w_0(\alpha_i) = \alpha_{i^*}$, where $\alpha_i$ is the
simple root corresponding to the vertex $i$. In the labeling in
\cite[Ch.~4]{Kac}, it is given by $i^* = \ell-i+1$ for type $A_\ell$,
$1^* = 5$, $2^* = 4$, $3^*=3$, $4^* = 2$, $5^* = 1$, $6^* = 6$ for
type $E_6$ respectively.
For type $D_\ell$ with odd $\ell$, it is given by
\begin{equation*}
    i^* =
    \begin{cases}
        \ell -1 & \text{if $i=\ell$}, \\
        \ell & \text{if $i=\ell-1$},\\
        i & \text{otherwise},
    \end{cases}
\end{equation*}
\begin{NB}
    For $D_\ell$ with odd $\ell$, it exchanges two `tails'. For $E_6$,
    it is a reflection at the center.
\end{NB}%
It is the identity for other types.

When $i^* = i$, we determine whether $\rho_i\cong\rho_i^*$ is given by
a symplectic or orthogonal form as follows: 
The trivial representation, assigned to $i=0$, is orthogonal.
For type $A_{\ell}$ with odd $\ell$ with $i= {(\ell+1)/2}$, it is
orthogonal. For other types, $\rho_{i_0}$ for the vertex $i_0$
adjacent to the vertex $i=0$ in the affine Dynkin diagram is the
representation given by the inclusion $\Gamma\subset\SU(2)$. Therefore
$\rho_{i_0}\cong \rho_{i_0}^*$ is symplectic. If $i$ is adjacent to
$i_0$ and $i^* = i$, then $\rho_i\cong\rho_i^*$ is orthogonal. If $j$
is adjacent to such an $i$ with $j^*=j$, then $\rho_j\cong\rho_j^*$ is
symplectic, and so on.

\subsection{ADHM description of dual instantons}\label{subsec:ADHM}

Let $(I,H)$ be the McKay quiver for $\Gamma$, namely $I$ is the set of
isomorphism classes of irreducible representations of $\Gamma$, and
$H$ is the set of arrows where we draw $a_{ij}$ arrows from $i$ to $j$
for $a_{ij} = \dim \Hom_\Gamma(\rho_i,\rho_j\otimes Q)$, where $Q$ is
the $2$-dimensional representation of $\Gamma$ given by the inclusion
$\Gamma\subset\SU(2)$. We also choose an orientation $\Omega$ of $H$,
which is a division $H = \Omega \sqcup\overline{\Omega}$, where
$\overline{h}$ is the arrow with the opposite direction to $h$. We define
$\varepsilon\colon H \to \{\pm 1\}$ as $\varepsilon(h) = 1$ if $h\in\Omega$,
$\varepsilon(h) = -1$ otherwise.

We take $\zeta = (\zeta_i) \in (\RR^3)^I$, the data for an ALE space
$X_\zeta$ asymptotic to $\RR^4/\Gamma$ at infinity. It sits in the
level $0$ hyperplane
$0 = \zeta\cdot\delta = \sum_i \zeta_i \dim \rho_i$, where $\delta$ is
the positive primitive imaginary root of the corresponding affine Lie
algebra.

Let us take an $\U(n)$ framed instanton $A$ on $X_\zeta$. We have the
corresponding ADHM description \cite{KN}. Namely we have $I$-graded
vector spaces $V=\bigoplus_{i\in I} V_i$, $W=\bigoplus_{i\in I} W_i$, and
linear maps $B_h\colon V_{\vout{h}}\to V_{\vin{h}}$ ($h\in H$),
$a_i\colon W_i\to V_i$, $b_i\colon V_i\to W_i$ ($i\in I$) satisfying
the hyperk\"ahler moment map equation. Here $\vout{h}$, $\vin{h}$ are
the outgoing and incoming vertices of an arrow $h$.
($B_h$ was denoted by $B_{i,j}$, and $a_i$, $b_i$ were denoted by
$i_k$, $j_k$ in \cite{KN}.)

Recall that the framing of $A$ is an approximate isomorphism of $A$
and a flat connection on $\RR^4/\Gamma$ at infinity. The flat
connection corresponds to a $\Gamma$-module
$\bigoplus_i W_i\otimes\rho_i$.

Recall that the description of \cite{KN} uses the tautological bundle
$\mathscr R$, which decomposes as
$\bigoplus_{i\in I} \mathscr R_i\otimes\rho_i^*$.
Reflection functors in \cite{Na-reflect} are understand as
isomorphisms between different descriptions of instanton moduli spaces
for different choices of $\mathscr R$. In particular, the reflection
functor $\mathcal F$ for the longest element $w_0$ composed with the
diagram automorphism $*$
\begin{equation*}
  * \circ \mathcal F = \mathcal F\circ * \colon
  \mathfrak M_\zeta^{\mathrm{reg}}(\bv,\bw)\to
  \mathfrak M_{w_0\zeta^*}^{\mathrm{reg}}(w_0\ast \bv^*,\bw^*)
\end{equation*}
corresponds to
$\mathscr R^* = \bigoplus_{i\in I}\mathscr R_i^*\otimes\rho_i$ as
shown in \cite[9(iii)]{Na-reflect}.
Here $\mathfrak M_\zeta^{\mathrm{reg}}(\bv,\bw)$ is the space of all
solutions of the hyperk\"ahler moment map equation modulo the action
of the group $\prod_{i\in I} \U(V_i)$, and $\bv = \dim V$,
$\bw = \dim W$ are dimension vectors of $V$, $W$. We also take
the open subset consisting of free $\prod_{i\in I} \U(V_i)$-orbits.
The reflection functor changes the dimension vector $\bv$ to $\bv'$ so
that $w_0(\bw - \mathbf C\bv) = \bw - \mathbf C \bv'$, and we denote
$\bv'$ by $w_0\ast \bv$. The parameter $\zeta$ is also changed as
$\zeta\to w_0\zeta$.
The diagram automorphism $*$ also changes dimension vectors, hence
$\mathfrak M_{w_0\zeta^*}^{\mathrm{reg}}(w_0\ast \bv^*,\bw^*)$ is the
space for the transformed dimension vectors and the parameter $\zeta$.
Note that $\mathcal F$ and $*$ do \emph{not} touch the vertex $0$
corresponding to the trivial representation, hence the component
$\bv_0$, $\bw_0$ are unchanged.

The dual instanton $A^*$ is given by first take the dual vector spaces
$V^* = \bigoplus_{i\in I} V_i^*$, $W^*=\bigoplus_{i\in I} W_i^*$, and then
replace linear maps as
\begin{equation*}
  -\varepsilon(h){}^t\! B_{\overline{h}}\colon V_{\vout{h}}^*\to V_{\vin{h}}^*,
  \quad
  - {}^t b_i \colon W_i^*\to V_i^*,
  \quad
  {}^t\! a_i\colon V_i^*\to W_i^*,
\end{equation*}
where the transpose ${}^t\xi\colon F^*\to E^*$ of a linear map
$\xi\colon E\to F$ is defined by $\langle \xi e, f\rangle = \langle e,
{}^t \xi f\rangle$ for $e\in E$, $f\in F^*$.
The parameter $\zeta$ for the hyperk\"ahler moment map equation is
replaced by $-\zeta$.
Hence this construction gives an involutive isomorphism
\begin{equation*}
  t\colon \mathfrak M_\zeta^{\mathrm{reg}}(\bv,\bw)\to
  \mathfrak M_{-\zeta}^{\mathrm{reg}}(\bv,\bw).
\end{equation*}

Composing two isomorphisms, we get
\begin{equation*}
  t\circ * \circ \mathcal F
  = *\circ\mathcal F\circ t\colon \mathfrak M_\zeta^{\mathrm{reg}}(\bv,\bw)
  \to \mathfrak M_\zeta^{\mathrm{reg}}(w_0\ast \bv^*,\bw^*),
\end{equation*}
where we have used $-w_0\zeta^* = \zeta$.

A little more precisely, the space $\mathfrak M_\zeta^{\mathrm{reg}}(\bv,\bw)$
depends on the $I$-graded vector space $W$, not only on its dimension
vector $\bw$. (On the other hand, it does not depend on $V$ as we take
the quotient by $\prod\U(V_i)$.) We replace
$W = \bigoplus_{i\in I} W_i$ by $\bigoplus_{i\in I} W_{i^*}^*$. Or
even better we should regard $W$ as a $\Gamma$ representation
$\bigoplus W_i\otimes\rho_i$, and it is replaced by its dual
representation
$\bigoplus W_i^*\otimes\rho_i^* = \bigoplus W_{i^*}^* \otimes\rho_i$.
It is more natural in view of \cite{KN}.

\subsection{$\SO/\grpSp$ instantons}

Let us take $\SO(n)$ or $\grpSp(n/2)$ framed instanton $A$ on
$X_\zeta$. We regard it as a $\U(n)$ framed instanton $A$ together
with an isomorphism $A\cong A^*$ compatible with the framing.
Since framed instantons have no nontrivial automorphisms, the
isomorphism $A\cong A^*$ is unique if it exists. Hence moduli spaces
of $\SO(n)/\grpSp(n/2)$ framed instantons are fixed point loci in
moduli spaces of $\U(n)$ instantons with respect to the involution
given by $A\mapsto A^*$. Here moduli spaces are constructed by a gauge
theoretic method as in \cite{MR1074476}.

When we say $A\cong A^*$ is compatible with the framing, we need to
fix an isomoprhism of the flat connection and its dual in advance. It
is given by a symplectic or nondegenerate symmetric form on the
representation $\bigoplus W_i\otimes\rho_i$, according to the
$\SO(n)$ or $\grpSp(n/2)$ instanton.

Let us describe the involution in terms of the ADHM description by
using the result explained in \subsecref{subsec:ADHM}. Since the gauge
theoretic construction and the ADHM description give isomorphic
hyperk\"haler manifolds \cite[\S8]{KN}, moduli spaces of
$\SO(n)/\grpSp(n/2)$ framed instantons are fixed point loci in
$\mathfrak M_\zeta^{\mathrm{reg}}(\bv,\bw)$ with respect to an
involution, which is more or less clear that it is given by
$t\circ * \circ\mathcal F$.

Note that $\bw^* = \bw$ as $\bigoplus W_i\otimes\rho_i$ is isomorphic
to its dual as above. The first Chern class as a $\U(n)$ instanton
vanishes, hence we have $w_0\ast \bv^* = \bv$ by \cite[\S9]{KN}.
But this is not enough to make $t\circ * \circ\mathcal F$ an
involution on $\mathfrak M_\zeta^{\mathrm{reg}}(\bv,\bw)$ as explained
in the last paragraph of \subsecref{subsec:ADHM}.
Namely we need to choose an isomorphism $W_i\cong W^*_{i^*}$. It is
given as follows.

Suppose $i=i^*$ (that is $\rho_i^*\cong \rho_i$). We have a symmetric
or symplectic form on $W_i\otimes\rho_i$ according to $\SO(n)$ or
$\grpSp(n/2)$ instantons. Since $\rho_i$ has a symmetric or symplectic
form according to the rule to determine the form for
$\rho_i\cong\rho_i^*$ is explained in \subsecref{subsec:diagram}, we
choose a symmetric or symplectic form on $W_i$. For example
$W_i\cong W_i^*$ is given by a symplectic form if
$\rho_i\cong\rho_i^*$ is symplectic and we consider $\SO(n)$
instantons.

When $i\neq i^*$, we have a symmetric or sympletic form on
$(W_i\otimes \rho_i) \oplus (W_{i^*}\otimes \rho_i^*)$.
Thus we have isomorphisms $f_i\colon W_i\to W_{i^*}^*$,
$f_{i^*}\colon W_{i^*}\to W_i^*$.
\begin{NB}
  They are defined by
  \begin{equation*}
    \langle f_i w_i, w_{i^*}\rangle = (w_i, w_{i^*}), \quad
    \langle f_{i^*} w_{i^*}, w_{i}\rangle = (w_{i^*}, w_{i}).
  \end{equation*}
  Then
  \begin{equation*}
    \langle {}^t\! f_{i^*} w_i, w_{i^*}\rangle
    = \langle w_{i^*}, {}^t\! f_{i^*} w_i\rangle
    = \langle f_{i^*} w_{i^*}, w_i\rangle
    = (w_{i^*},w_i) = \pm (w_i, w_{i^*})
    = \pm \langle f_i w_i, w_{i^*}\rangle.
  \end{equation*}
\end{NB}%
Then we have ${}^t\! f_{i^*} = \pm f_i$, where we choose the sign
according to $\SO(n)$ or $\grpSp(n/2)$ instantons.

Now $t\circ * \circ\mathcal F$ is defined as an involution on
$\mathfrak M_\zeta^{\mathrm{reg}}(\bv,\bw)$.
\begin{Theorem}
  A moduli space of $\SO(n)$ or $\grpSp(n/2)$ framed instantons on
  $X_\zeta$ is isomorphic to the fixed point locus
  $\mathfrak M_\zeta^{\mathrm{reg}}(\bv,\bw)^{t\circ *\circ\mathcal
    F}$ as a hyperk\"ahler manifold.
\end{Theorem}

Dimension vectors $\bv$, $\bw$ are determined by the second Chern
class and the framing considered as $\U(n)$ instantons.

\subsection{Extension to partial compactification}

Recall that $\mathfrak M_\zeta^{\mathrm{reg}}(\bv,\bw)$ is defined as
the open subset consisting of free $\prod \U(V_i)$-orbits. We have a
larger space $\mathfrak M_\zeta(\bv,\bw)$ by dropping the freeness
condition. It is Uhlenbeck's partial compactification of the moduli
space of framed instantons on $X_\zeta$.

The isomorphism $t\circ * \circ\mathcal F$ extends to
$\mathfrak M_\zeta(\bv,\bw)$ as a homeomorphism. This is because the
reflection functor is defined on the larger space, and $t$, $*$
clearly extend. It is also clear in the gauge theoretic construction
of $\mathfrak M_\zeta(\bv,\bw)$.
If we use algebro-geometric construction of
$\mathfrak M_\zeta(\bv,\bw)$ via geometric invariant theory, we can
endow $\mathfrak M_\zeta(\bv,\bw)$ with a structure of a
quasiprojective variety. Then the extension of
$t\circ * \circ\mathcal F$ is an involution on a variety. The fixed
point loci are also quasiprojective varieties.

Note that $\mathfrak M_\zeta^{\mathrm{reg}}(\bv,\bw)$ has another
partial compactification as a moduli space of framed torsion free
sheaves on $X_\zeta$. The corresponding ADHM description is given as
follows. (See \cite{Na-ADHM} for detail.) We take an algebro-geometric
description of $\mathfrak M_\zeta(\bv,\bw)$. We
decompose the parameter $\zeta$ to complex and real parts $\zeta_\CC$,
$\zeta_\RR$, impose the complex moment map equation involving only
$\zeta_\CC$. We have the group action of $\prod \GL(V_i)$ on the
solution space. We take the quotient of the $\zeta_\RR$-semistable
locus by the S-equivalene relation. Thus
\begin{equation*}
  \mathfrak M_\zeta(\bv,\bw) \cong H^{\mathrm{ss}}_{(\zeta_\RR,\zeta_\CC)}/\!\sim,
\end{equation*}
where $H^{\mathrm{ss}}_{(\zeta_\RR,\zeta_\CC)}$ is the open subset of
$\zeta_\RR$-semistable points in the solution space of the complex
moment map equation. See \cite[Prop.~2.11]{Na-reflect}.

Recall that $\zeta_\RR$ lives on the level 0 hyperplane
$\zeta_\RR\cdot \delta = 0$. Then we take $\zeta_\RR'$ near
$\zeta_\RR$ with $\zeta_\RR'\cdot\delta < 0$. Then we define
$\mathfrak M_{(\zeta_\CC,\zeta_\RR')}(\bv,\bw)$ as the quotient of the
$\zeta_\RR'$-semistable (equivalently $\zeta_\RR'$-stable) locus by
the action of $\prod\GL(V_i)$.
Then $\zeta_\RR'$-semistability implies $\zeta_\RR$-semistability by
our choice, that is $H^{\mathrm{ss}}_{(\zeta_\RR',\zeta_\CC)}\subset
H^{\mathrm{ss}}_{(\zeta_\RR,\zeta_\CC)}$. It induces a morphism
\begin{equation*}
  \mathfrak M_{(\zeta_\CC,\zeta_\RR')}(\bv,\bw) \to
  \mathfrak M_\zeta(\bv,\bw), 
\end{equation*}
which is an isomorphism on
$\mathfrak M_\zeta^{\mathrm{reg}}(\bv,\bw)$.  Therefore
$\mathfrak M_{(\zeta_\CC,\zeta_\RR')}(\bv,\bw)$ is another partial
compactification of $\mathfrak M_\zeta^{\mathrm{reg}}(\bv,\bw)$.

Note that $-w_0 (\zeta_\RR')^* \neq \zeta_\RR'$ as they live in the
opposite side of the level 0 hyperplane $\zeta_\RR\cdot\delta = 0$.
Therefore $t\circ * \circ\mathcal F$ does not define an involution
on $\mathfrak M_{(\zeta_\CC,\zeta_\RR')}(\bv,\bw)$.

\bibliographystyle{myamsalpha}
\bibliography{nakajima,mybib,coulomb}

\def\cprime{$'$} \def\cprime{$'$} \def\cprime{$'$} \def\cprime{$'$}
  \def\cprime{$'$}
  \providecommand{\noopsort}[1]{}\def\cftil#1{\ifmmode\setbox7\hbox{$\accent"5E#1$}\else
  \setbox7\hbox{\accent"5E#1}\penalty 10000\relax\fi\raise 1\ht7
  \hbox{\lower1.15ex\hbox to 1\wd7{\hss\accent"7E\hss}}\penalty 10000
  \hskip-1\wd7\penalty 10000\box7}
\providecommand{\bysame}{\leavevmode\hbox to3em{\hrulefill}\thinspace}
\providecommand{\MR}{\relax\ifhmode\unskip\space\fi MR }
% \MRhref is called by the amsart/book/proc definition of \MR.
\providecommand{\MRhref}[2]{%
  \href{http://www.ams.org/mathscinet-getitem?mr=#1}{#2}
}
\providecommand{\href}[2]{#2}
\begin{thebibliography}{AHDM78}

\bibitem[ADHM78]{ADHM}
M.~F. Atiyah, N.~J. Hitchin, V.~G. Drinfel{\cprime}d, and Y.~I. Manin,
  \emph{Construction of instantons}, Phys. Lett. A \textbf{65} (1978), no.~3,
  185--187. \MR{598562 (82g:81049)}

\bibitem[Buc86]{MR857374}
N.~P. Buchdahl, \emph{Instantons on {${\bf C}{\rm P}_2$}}, J. Differential
  Geom. \textbf{24} (1986), no.~1, 19--52. \MR{857374}

\bibitem[Cho16]{MR3508922}
J.~Choy, \emph{Moduli spaces of framed symplectic and orthogonal bundles on
  {$\Bbb{P}^2$} and the {$K$}-theoretic {N}ekrasov partition functions}, J.
  Geom. Phys. \textbf{106} (2016), 284--304. \MR{3508922}

\bibitem[DK90]{MR1079726}
S.~K. Donaldson and P.~B. Kronheimer, \emph{The geometry of four-manifolds},
  Oxford Mathematical Monographs, The Clarendon Press, Oxford University Press,
  New York, 1990, Oxford Science Publications. \MR{1079726 (92a:57036)}

\bibitem[Kac90]{Kac}
V.~G. Kac, \emph{Infinite-dimensional {L}ie algebras}, third ed., Cambridge
  University Press, Cambridge, 1990. \MR{MR1104219 (92k:17038)}

\bibitem[KN90]{KN}
P.~B. Kronheimer and H.~Nakajima, \emph{Yang-{M}ills instantons on {ALE}
  gravitational instantons}, Math. Ann. \textbf{288} (1990), no.~2, 263--307.
  \MR{MR1075769 (92e:58038)}

\bibitem[{Li}18]{QuiverSym}
Y.~{Li}, \emph{{Quiver varieties and symmetric pairs}}, ArXiv e-prints (2018),
  \href{http://arxiv.org/abs/1801.06071}{{\ttfamily arXiv:1801.06071
  [math.RT]}}.

\bibitem[Nak90]{MR1074476}
H.~Nakajima, \emph{Moduli spaces of anti-self-dual connections on {ALE}
  gravitational instantons}, Invent. Math. \textbf{102} (1990), no.~2,
  267--303. \MR{MR1074476 (92a:58024)}

\bibitem[Nak03]{Na-reflect}
\bysame, \emph{Reflection functors for quiver varieties and {W}eyl group
  actions}, Math. Ann. \textbf{327} (2003), no.~4, 671--721. \MR{MR2023313
  (2004k:16036)}

\bibitem[Nak07]{Na-ADHM}
\bysame, \emph{Sheaves on {$ALE$} spaces and quiver varieties}, Mosc. Math. J.
  \textbf{7} (2007), no.~4, 699--722. \MR{MR2372210}

\bibitem[Nak16]{2015arXiv150303676N}
\bysame, \emph{{Towards a mathematical definition of Coulomb branches of
  $3$-dimensional $\mathcal N=4$ gauge theories, I}}, Adv. Theor. Math. Phys.
  \textbf{20} (2016), no.~3, 595--669,
  \href{http://arxiv.org/abs/1503.03676}{{\ttfamily arXiv:1503.03676
  [math-ph]}}.

\end{thebibliography}

\end{document}